\documentclass[12pt,reqno,oneside]{amsart}
\usepackage{a4}
\usepackage{amssymb}
\usepackage{lcy}
\usepackage{verbatim}
\usepackage[russian]{babel}

\newtheorem{prp}{}[section]
\newtheorem{que}{Вопрос}

\def\iat{\operatorname{IA}}

\def\gl#1{\mathop{\rm GL}(#1)}

\def\aut#1{\mathop{\rm Aut}(#1)}
\def\out#1{\mathop{\rm Out}(#1)}

\def\id{\mathop{\rm id}}
\def\ia#1{\mathop{\rm IA}(#1)}
\def\inn#1{\mathop{\rm Inn}(#1)}

\def\inv{^{-1}}
\def\str#1{\langle#1\rangle}

\def\f{\varphi}
\def\a{\alpha}

\def\s{\sigma}

\def\Theo{\mathop{\rm Th}}
\def\Q{\mathbf Q}

\def\cB{{\mathcal B}}

\def\cD{{\mathcal D}}

\def\cX{\mathcal X}
\def\cY{\mathcal Y}

\def\cM{\mathcal M}
\def\cL{\mathcal L}
\def\cN{\mathcal N}

\def\N{\mathbf N}
\def\Z{\mathbf Z}
\def\Q{\mathbf Q}
\def\R{\mathbf R}

\renewcommand{\le}{\leqslant}
\renewcommand{\ge}{\geqslant}

\def\sym#1{\mathop{\rm Sym}(#1)}
\def\rank{\mathop{\rm rank}\,}

\def\hatin{\widehat{\phantom a}}
\def\restr{\!\upharpoonright\!}

\def\pw#1{{}_{\{#1\}}}
\def\sw#1{{}_{(#1)}}

\begin{document}

УДК 512.543.25, 510.67


\bigskip

\title[Группы автоморфизмов отн. свободных групп]{Группы автоморфизмов относительно свободных групп бесконечного ранга}
\author{В.~А.~Толстых}
\maketitle

Настоящий обзор посвящен недавним результатам о группах
автоморфизмов относительно свободных групп бесконечного
ранга. Краткие исторические замечания и ссылки
на литературу, которые мы приводим при обсуждении
более общих или смежных тем, не претендуют
на полноту.

\section{Топологические аспекты}

Пусть $\Omega$ --- бесконечное множество. Множество
$\Omega^\Omega$ всех отображений $\Omega$ в себя
превращается в топологическое пространство,
если в качестве предбазы топологии, {\it функциональной
топологии} на $\Omega^\Omega,$ рассмотреть семейства
функций вида $\{\s \in \Omega^\Omega : \s x =y\},$
где $x,y$ пробегают $\Omega$ \cite{BSh},\cite[глава 4]{Ho}.  Легко видеть, что
операция композиции на $\Omega^\Omega$ непрерывна в
функциональной топологии, и что сужая функциональную
топологию на любую группу перестановок $\Gamma$
множества $\Omega,$ мы превращаем $\Gamma$ в
топологическую группу. Мы будем использовать
стандартные обозначения теории бесконечных групп
перестановок: символ $\Gamma\pw U,$ где
$U$ --- подмножество $\Omega,$ обозначает
подгруппу всех элементов $\Gamma,$
стабилизирующих $U$ поточечно, а $\Gamma\sw U$ ---
подгруппу всех элементов $\Gamma,$
стабилизирующих $U$ как множество.

Естественен вопрос о том, можно ли восстановить
топологию на $\Gamma,$ используя исключительно
теоретико-групповые средства. Например, если
$\Gamma$ --- группа автоморфизмов структуры $\cM,$
определенной на $\Omega,$ то восстановление
функциональной топологии на $\Gamma$ может
служить первым этапом в восстановлении
$\cM$ в $\Gamma$.

В работе \cite{DiNeuTho} Дж.~Диксона, П.~Ноймана,
С.~Томаса и работе \cite{MN} П.~Ноймана и
Д.~Макферсона был получен положительный ответ для
случая $\Gamma=\sym \Omega.$
В работе \cite{DiNeuTho} рассматривался случай, когда $\Omega$ ---
счетное множество
(в этом случае функциональная топология
наиболее полезна, см. ниже), а в работе \cite{MN} --- общий случай.
В работе \cite{DiNeuTho}
было установлено, что если $\Omega$ --- счетное
множество, то подгруппа
$G$ группы $\sym \Omega$ открыта в том и только
в том случае, если индекс $G$ в группе $\sym \Omega$
меньше кардинала $2^{\aleph_0}.$ Более общий
результат из работы \cite{MN} звучит так: если $\Omega$ ---
произвольное бесконечное множество, то
подгруппа $G$ группы $\sym \Omega$ является
открытой в том и только в том случае, если
$\sym \Omega = \text{FS}(\Omega) G,$
где $\text{FS}(\Omega)$ --- группа
всех финитарных перестановок множества $\Omega.$
Перестановка $\s \in \sym \Omega$ называется
{\it финитарной}, если ее носитель
$
\text{supp}(\s) = \{x \in \Omega : \s x \ne x\}
$
конечен.

Пусть $G$ --- бесконечнопорожденная группа. {\it
Конфинальностью} $G$ называется мощность наименьшей
цепочки собственных подгрупп, объединение которой есть
$G.$ Еще один важный результат из \cite{MN} устанавливает,
что конфинальность бесконечной симметрической группы
$\sym \Omega$ строго больше $|\Omega|.$

Отличительным свойством группы $\sym \Omega$ над
счетно-бесконечным множеством $\Omega,$ наделенной
функциональной топологией,  является то, что она
является {\it польской группой} \cite[раздел 9]{Kech}.
Это означает, что существует метрика на $\sym \Omega,$
совместимая с функциональной топологией, и по
отношению к этой метрике $\sym \Omega$ является полным
и сепарабельным метрическим пространством. К примеру,
если $\Omega = (a_n : n \in \N),$ то совместимая
метрика $d,$ обладающая необходимыми свойствами, может
быть задана так: $d(\s,\pi) =0,$ если $\s=\pi$ и
если $\s \ne \pi,$ то
$d(\s,\pi)=2^{-n},$ где $n$ --- наименьшее натуральное
число, такое, что $\s(a_n) \ne \pi(a_n)$ или
$\s\inv(a_n) \ne \pi\inv(a_n).$ Можно, далее,
показать, что группы автоморфизмов всех
структур первого порядка, определенных на
$\Omega,$ и только они являются замкнутыми
подгруппами $\sym \Omega$ \cite[глава 4]{Ho}. Замкнутая подгруппа
польской группы снова является польской, и потому
группы автоморфизмов всех счетно-бесконечных
структур являются польскими группами.

Упомянутый выше результат из работы \cite{DiNeuTho} стимулировал поиск
примеров счетных структур первого порядка $\cM,$
таких, что открытые подгруппы их групп автоморфизмов
$\aut{\cM}$ описываются тем же условием, что и
открытые подгруппы симметрической группы на счетном
бесконечном множестве, т.е. в точности как подгруппы,
имеющие "малый", меньше $2^{\aleph_0},$ индекс.  Про
всякую такую счетную структуру говорят, что она имеет
{\it свойство малого индекса.} Среди примеров структур,
имеющих свойство малого индекса, можно указать
векторные пространства счетно-бесконечной размерности
над не более, чем счетными телами, структуру
$\str{\Q;<},$ счетные безатомные булевы алгебры и т.д.

В.~Ходжес, И.~Ходкинсон, Д.~Ласкар и С.~Шелах показали
в работе \cite{Hodges_et_al}, что счетные $\omega$-категоричные
и $\omega$-стабильные структуры, а также так называемый
случайный граф (предел Фрессе класса всех конечных
графов, см. \cite[глава 7]{Ho}) имеют свойство малого индекса,
а их группы автоморфизмов --- несчетную конфинальность.
Р.~Брайан и Д.~Эванс установили в \cite{BrEv}, что ключевые
методы из статьи  \cite{Hodges_et_al} могут
быть применены ко многим относительно свободным
группам счетно-бесконечного ранга, весьма далеким
во многих отношениях от структур, изучавшихся
в \cite{Hodges_et_al}. А именно, Р.~Брайан и Д.~Эванс
установили, что при выполнении относительно
свободной группой $G$ счетно-бесконечного
ранга некоторого естественного свойства, названного ими свойством
базисной конфинальности, группа $\aut G$ относительно
свободной группы $G$ счетного ранга {богата
так называемыми генерическими автоморфизмами}.

Определение генерических автоморфизмов относительно
свободных групп, данное Р.~Брайаном и Д.~Эвансом, следует
идее из \cite{Hodges_et_al}. Пусть $G$ ---
относительно свободная группа счетно-бесконечного
ранга, $\Gamma$ обозначает группу $\aut G,$ а $\cB(G)$
--- семейство всех свободных множителей группы
$G,$ имеющих конечный ранг. Говорим, что набор
$(\gamma_1,\ldots,\gamma_n)$ элементов $\Gamma^n$
является {\it $\cB(G)$-генерическим} или просто {\it
генерическим}, если выполнены следующие два условия:

(1) для каждого свободного множителя $A \in \cB(G)$
подгруппы $\Gamma_{(B)},$ где $A \subseteq B \in \cB(G)$ и $\gamma_i
B=B$ для всех $i=1,\ldots,n,$
формируют базу открытых окрестностей единицы;

(2) если $A \in \cB(G),$ $\gamma_i A =A$ для всех
$i=1,\ldots,n,$ $B$ --- свободный множитель из $\cB(G),$
содержащий $A,$ $\beta_1,\ldots,\beta_n\in \aut B$ и $\beta_i\restr A
= \gamma_i\restr A$ для всех $i=1,\ldots,n,$ то найдется
автоморфизм $\alpha \in \Gamma_{(A)},$ такой, что
$\gamma^\alpha_i \restr B =\beta_i$ для всех
$i=1,\ldots,n$ (символ $\restr$ указывает
на ограничения отображений, а $\gamma^\alpha=\alpha
\gamma \alpha\inv$).

Функциональная топология на $\Gamma$ естественным образом задает
топологию на декартовых степенях
$\Gamma^n.$ Группа $\Gamma=\aut G$ {\it богата
генерическими автоморфизмами} \cite{BrEv}, если для
каждого натурального $n$ множество всех наборов из
$\Gamma^n,$ которые являются генерическими --- котощее
(множество второй категории) в $\Gamma^n.$

Свойство базисной конфинальности, которое как мы
указали выше, является достаточным условием для того,
чтобы группа автоморфизмов относительно
свободной группы $G$ счетно-бесконечного ранга была богата генерическими
автоморфизмами формулируется так: пусть $\{y_n : n \in \N\}$
--- какой-нибудь базис группы $G$; тогда
для каждого $\alpha\in
\Gamma$ и каждого натурального число $n$
найдутся натуральное число $r\ge n$  и автоморфизм $\beta \in {\rm
Aut}(\langle y_1,\ldots, y_r\rangle)$, такие, что
$\beta y_i=\alpha y_i$ для всех $i=1,\ldots,n$. В работе
\cite{BrEv} показано, что если $G$ --- абсолютная
свободная группа или если $\aut G$ --- ручная группа,
то $G$ имеет свойство базисной конфинальности.

Пусть $\mathfrak V$ --- многообразие
групп, а $F$ --- абсолютно свободная группа. Рассмотрим относительно свободную группу
$F/\mathfrak V(F)$ многообразия $\mathfrak V,$ где $\mathfrak
V(F)$ --- вербальная подгруппа $F,$ определяющая
многообразие $\mathfrak V.$ Тогда группа автоморфизмов
$\aut{F/\mathfrak V(F)}$ называется {\it ручной},
если  гомоморфизм $\aut F \to \aut{F/\mathfrak
V(F)},$ индуцированный естественным гомоморфизмом $F \to
F/\mathfrak V(F),$ является сюръективным.  К примеру,
согласно результату Р.~Брайана и О.~Мацедонской из \cite{BrMa},
для всякого многообразия $\mathfrak V,$
свободные группы которого нильпотентны,
группы автоморфизмов {\it бесконечнопорожденных}
групп ${F/\mathfrak V(F)}$ --- ручные
\cite{BrMa}; список известных ручных
групп автоморфизмов бесконечнопорожденных относительно свободных
приводится в \cite[теорема 1.6]{BrEv}. Существуют
примеры групп автоморфизмов относительно свободных
групп счетно-бесконечного ранга, не являющихся ручными
\cite{BrEv}.

Р.~Брайан и Д.~Эванс указывают \cite{BrEv}, что если
группа $\aut G$ счетной относительно свободной группы
$G$ богата генерическими автоморфизмами, то, рассуждая
практически так же как и в \cite{Hodges_et_al}, можно
вывести из этого обстоятельства, что группа $\aut G$
имеет свойство малого индекса и несчетную конфинальность.

Через $2^n$ мы обозначаем множество всех
отображений ординала $n \in \omega$
в $2=\{0,1\},$ $2^{<\omega}$ --- это
множество $\bigcup_{n \in \omega} 2^n,$
а $2^\omega$ --- множество всех отображений
$\omega$ в $2.$ Если $s \in 2^n,$ то
$s\widehat~0$ (соотв. $s\widehat~1$) ---
это продолжение $s$ на $n+1,$ принимающее
значение $0$ (соотв. $1$) на $n.$

Пусть $\cM$ --- счетная $\omega$-категоричная и
$\omega$-стабильная структура. Идея доказательств
результатов из \cite{Hodges_et_al}, упомянутых выше,
состоит в построении индукцией по элементам $s$
бинарного дерева $2^{<\omega}$ генерических
автоморфизмов $\gamma_{s}, g_{s\widehat~0},
g_{s\widehat~1}$ таких, что:
\begin{itemize}
\item[(i)] для
каждого $\sigma \in 2^\omega$ последовательность
$(\gamma_{\sigma \upharpoonright n})$ должна быть последовательностью
Коши;

\item[(ii)] если $\gamma_\pi$ обозначает
$\lim \gamma_{\pi \upharpoonright n},$ где $\pi \in 2^\omega,$
то $\gamma_\sigma \ne \gamma_\pi$
для всех $\sigma, \pi \in 2^\omega$ и

\item[(iii)] из условий $\sigma \ne \tau,$
$\sigma \upharpoonright n =s\widehat~0,$
$\tau \upharpoonright n =s\widehat~1,$ где
$\sigma,\tau \in 2^\omega,$ а $n \in \omega,$ должно
следовать, что
\begin{equation} \label{0-1Cond}
(g_{s\widehat~0})^{\gamma_\sigma \gamma_\tau\inv} =g_{s\widehat~1}.
\end{equation}

\end{itemize}
При этом одно из условий на автоморфизмы $g_s,$ которое
может быть названо 0-1 условием, варьируется
в зависимости от того, какой цели мы хотим добиться.
К примеру, для того, чтобы доказать (от противного),
что подгруппа $H$ группы $\aut \cM,$ имеющая малый
индекс, открыта, $g_{s\widehat~0}$ выбирают
в $H,$ а $g_{s\widehat~1}$ --- вне $H.$ Тогда
из (\theequation) следует, что $H$ имеет
индекс $2^{\aleph_0},$ противоречие.

И.~В.~Чирков \cite{Chir}, следуя идеям
Р.~Брайана и Д.~Эванса из \cite{BrEv}, показал, что
всякая свободная алгебра Ли счетно-бесконечного ранга тоже имеет
свойство малого индекса, а ее группа автоморфизмов
--- несчетную конфинальность.

Пусть $\cX$ --- базис относительно свободной
группы $F$ бесконечного ранга. Группа
$\text{Aut}_{\text{fin},\cX}(F)$ состоит из всех
автоморфизмов $F,$ действительно перемещающих
не более конечного числа элементов $\cX.$

\begin{que}
Пусть $F$ --- свободная {\em(}свободная
абелева, свободная нильпотентная, свободная
разрешимая{\em)} группа бесконечного ранга,
а $\cX$ --- это некоторый базис $F.$ Верно ли,
что подгруппа $H$ группы $\aut F$ открыта
в том и только в том случае, если
$\aut F = H \text{\rm Aut}_{\text{\rm fin},\cX}(F)?$
\end{que}

Иначе говоря, мы требуем, чтобы открытые подгруппы
$\aut F$ допускали описание в духе приведенного выше
результата из \cite{MN} об открытых подгруппах
произвольных бесконечных симметрических групп.

Применяя топологические методы, Р.~Брайан
и В.~А.~Романьков \cite{BrRom} получили
ряд интересных общих результатов о группах
автоморфизмов относительно свободных алгебр (т.е. свободных
объектов из многообразий алгебр --- структур первого
порядка, в язык которых входят только функциональные
символы).

Первый круг результатов из \cite{BrRom} относится к
свободным подгруппам группах автоморфизмов относительно
свободных алгебр бесконечного ранга. Обобщая известные
результаты о симметрических группах, Р.~Брайан и
В.~А.~Романьков показывают, что если $F$ --- относительно
свободная алгебра бесконечного ранга $\varkappa,$ то
группа $\Gamma=\aut F$ имеет {\it всюду плотные} свободные
подгруппы любой мощности $\nu,$ где $|F| \le \nu
\le 2^\varkappa.$ Еще один результат о свободных
подгруппах $\aut F$ говорит, неформально, о том, что
"большинство"\ конечнопорожденных подгрупп $\aut F$ ---
это свободные группы. А именно,  множество $\Phi_n$
всех наборов $(\varphi_1,\ldots,\varphi_n)$ декартовой
степени $\Gamma^n,$ таких, что подгруппа
$\langle\varphi_1,\ldots,\varphi_n\rangle,$
порожденная автоморфизмами
$\varphi_1,\ldots,\varphi_n,$ --- свободная группа
ранга $n,$ является котощим в $\Gamma^n.$

Пусть $\pi$ --- автоморфизм относительно свободной
алгебры $F$ бесконечного ранга, действующий
на некотором базисе как перестановка с бесконечными
циклами и без неподвижных элементов. Тогда,
как показано в \cite{BrRom}, нормальное замыкание
$\pi$ в $\aut F$ является всюду плотной подгруппой
$\aut F.$ В некоторых случаях указанное
нормальное замыкание может, конечно, совпадать с $\aut F.$
Это имеет место, например, если $F$ --- алгебра в пустом языке,
если $F$ --- это векторное пространство
над телом \cite{To_GL(V)}, если $F$ ---
свободная нильпотентная группа бесконечного ранга
\cite{To_BP} и т.п. Ниже мы еще вернемся к обсуждению этой темы.

Следующий результат имеет ряд интересных
теоретико-групповых следствий \cite{BrRom}.

\begin{prp} \label{BrRomGenSet}
Пусть $F$ --- относительно свободная алгебра
бесконечного ранга, $\cX$ --- некоторый базис
$F,$ $\Gamma=\aut F$ и $\sym \cX$ --- множество
всех автоморфизмов $F,$ действующих на $\cX$
как перестановки. Тогда
$$
\aut F = \str{\pi \Gamma_{(U)} \pi\inv: \pi \in \sym \cX},
$$
где $U$ --- произвольное конечное подмножество
$F.$ В частности, $\Gamma$ не имеет собственных
нормальных открытых подгрупп.
\end{prp}

Пусть теперь $F$ --- относительно свободная
алгебра бесконечного ранга, имеющая свойство малого индекса. Если
бы группа $\aut F$ имела собственную нормальную подгруппу
индекса $< 2^{\aleph_0},$ то такая подгруппа
была бы открытой, что невозможно по \ref{BrRomGenSet}.
Можно показать, что всякая абелева группа
имеет не более чем счетный гомоморфный
образ. Следовательно, группа $\aut F$ обязана
совпадать со своим коммутантом \cite{BrRom}, ибо
иначе $\aut F$ можно было бы гомоморфно отобразить
на счетную нетривиальную абелеву группу, что,
как мы видели, невозможно.

Пусть $L$ --- счетный язык, а $\frak V_1 \subseteq \frak V_2$
--- некоторые многообразия алгебр в языке $L.$ Пусть, далее,
$G$ --- свободная алгебра из $\frak V_1,$ а $F$ --- свободная
алгебра из $\frak V_2,$ причем обе алгебры $F,G$ имеют
счетно-бесконечный ранг. Выбирая в $F$ и $G$ некоторые
базисы, строим сюръективный гомоморфизм $F \to G,$ который, очевидно,
индуцирует гомоморфизм $\tau : \aut F \to \aut G$ групп
автоморфизмов. Пусть $T$ --- образ $\tau$; автоморфизмы
из группы $T$ естественно называть ручными (относительно
многообразия $\frak V_2$). Ясно, что
$\tau$ зависит от выбора базисов в $F,G,$ но после перемены
базисов образ индуцированного гомоморфизма будет
сопряжен с $T.$ Нормальное замыкание $T$
является всюду плотной подгруппой $\aut G,$ ибо сужение
$\tau$ на подгруппу $\aut F,$ состоящую из
автоморфизмов, действующих как перестановки на
выбранном базисе $F,$ является изоморфизмом.

Используя то обстоятельство, что группы $\aut F$ и
$\aut G$ --- польские, Р.~Брайан и В.~А.~Романьков описывают в
\cite{BrRom} следующую альтернативу для подгруппы $T$ ручных
автоморфизмов: либо $T=\aut G,$ либо $T$ --- тощее
подмножество $\aut G$ (множество первой категории),
причем в последнем случае индекс $T$ в $\aut G$ равен
$2^{\aleph_0}.$ Схема доказательства такова:
замечаем, что $\tau$ является непрерывным отображением,
и поэтому образ польской группы $\aut F$ в $\aut G$
--- это, по крайней мере, множество со свойством
Бэра. Согласно одному результату К.~Куратовского,
всякая подгруппа со свойством Бэра в польской группе
либо открыта, либо является тощим множеством. В первом случае
получаем по \ref{BrRomGenSet}, что $T=\aut G,$
ибо $T$ содержит подгруппу автоморфизмов,
действующих перестановками на выбранном
выше базисе $G.$ Во втором случае применяем
результат из \cite{Hodges_et_al}, гласящий,
что тощая подгруппа польской группы имеет
индекс $2^{\aleph_0}.$

\section{Порождающие множества}

Шутливое замечание из \cite[стр. 44]{LSch} гласит, что описание
порождающих группы автоморфизмов бесконечнопорожденной
свободной группы "представляет трудности лишь в
смысле обозначений". Тем не менее, о порождающих
множествах групп автоморфизмов бесконечнопорожденных
относительно свободных групп известно очень мало, что,
конечно, является серьезным препятствием при изучении этих групп.

Пожалуй, самым общим известным фактом является
процитированный выше результат из статьи \cite{BrRom}
Р.~Брайана и В.~А.~Романькова о порождении группы
автоморфизмов $\Gamma=\aut F$ относительно свободной
алгебры $F$ бесконечного ранга объединением подгрупп,
сопряженных подгруппе $\Gamma_{(U)}$ посредством
"перестановочных"\ автоморфизмов $F,$ где $U$
--- произвольное конечное подмножество $F$ (см.
\ref{BrRomGenSet}).

Для описания порождающих группы автоморфизмов $\Gamma$
бесконечнопорожденной свободной абелевой группы $A$
простой экспоненты $p$ можно применить результат А.~Розенберга
\cite{Rosen} о бесконечномерных общих линейных группах
над телами:  $\Gamma$ порождается всеми элементами
{\it класса два}, т.е. такими элементами $\s \in
\Gamma,$ что $(\s-\id)^2=0.$

Р.~Г.~Суон (см. \cite{BurnsPi}) нашел единообразно
описываемые множества порождающих как для групп
автоморфизмов бесконечнопорожденных свободных абелевых
групп, так и для групп автоморфизмов
бесконечнопорожденных свободных абелевых групп
произвольной экспоненты $m \ge 2.$
Пусть $A$ --- свободная абелева группа (свободная абелева группа экспоненты
$m$) и $\cX$ --- это произвольный базис $A.$ Тогда автоморфизм
$\theta$ группы $A$ называется
{\it $\cX$-блочно-унитреугольным}, если существует
половинное подмножество $\cY$ базиса $\cX,$ такое, что
$\theta$ оставляет на месте все элементы $\cY$ и
для каждого  $z \in \cX \setminus \cY$
$$
\theta z \equiv z\,(\operatorname{mod}\, \str{\cY}).
$$
(в частности, `матрица' $\theta$
в базисе $\cY \cup (\cX \setminus \cY)$
является блочно-унитреугольной;
говорят, что подмножество $M$ данного бесконечного
множества является {\it половинным} --- "moiety"\ в
англ. литературе, --- если мощность $M$ равна
мощности дополнения $M$ до этого множества.)

\begin{prp}[Р.~Г.~Суон] \label{Swan}
Группа $\aut A$ порождается всеми
$\cX$-блочно-уни\-тре\-уголь\-ными
автоморфизмами, и, более того,
ширина $\aut A$ относительно всех
этих автоморфизмов не превышает $22.$
\end{prp}

Напомним, что {\it ширина} данной группы $G$
относительно порождающего множества $S$
--- это наименьшее натуральное
число $k,$ такое, что каждый элемент группы
$G$ может быть записан в виде произведения
не более $k$ элементов из $S \cup S\inv,$
или $\infty,$ в противном
случае.

Недавний неожиданный результат Дж.~Бергмана
\cite{Berg} о бесконечных симметрических группах, а
также предложенные им проблемы, стимулировали
появление ряда работ о свойствах порождающих множеств групп
автоморфизмов различных структур. Пусть $\Omega$
--- произвольное бесконечное множество. Дж.~Бергман показал
в \cite{Berg}, что ширина группы $\sym \Omega$
относительно {\it любого} множества порождающих
конечна. Мы называем группу $G$ {\it группой конечной
ширины}, если ее ширина относительно любого
порождающего множества конечна (заметим, что используются
также термины как "свойство Бергмана",
"группы ограниченного диаметра Кэли"\ и т.п.).

По-видимому, первый пример бесконечной группы
конечной ширины был найден С.~Шелахом
в работе \cite{Sh} 1980 года.
Эта работа содержит пример несчетной группы
$G,$ которая имеет ширину не более чем 240
по отношению к каждому своему порождающему
множеству.

В первой версии препринта \cite{Berg_arxive}, предшествующем статье \cite{Berg},
Дж.~Бергман сформулировал несколько вопросов о том,
являются ли группы автоморфизмов различных
классических структур группами конечной
ширины. В частности. он предложил проанализировать ситуацию
для группы автоморфизмов множества вещественных
чисел $\R$ как {борелевского пространства,}
для групп автоморфизмов однородных булевых пространств,
для бесконечномерных общих линейных групп, для
групп автоморфизмов свободных групп бесконечного
ранга, а также для ряда других групп автоморфизмов.
Более того, Дж.~Бергманом была высказана общая гипотеза
о том, что новые примеры групп конечной ширины
могут быть найдены "среди групп автоморфизмов
структур, которые могут быть `собраны' из бесконечного
числа копий самих себя".

Вскоре после появления препринта Дж.~Бергмана М.~Дросте
и Р.~Гёбель \cite{DrGo} ответили положительно на вопрос
Дж.~Бергмана о группе автоморфизмов $\R$ как борелевского
пространства, а автор получил \cite{To_GL(V)}, тоже положительный,
ответ для случая бесконечномерных общих линейных
групп.  Другие примеры групп автоморфизмов, являющихся
группами конечной ширины, могут быть найдены в работах \cite{DrHo,KechRosen}.

{\it Степенью} $X^k$ порождающего множества $X$ данной
группы $G$ будем называть множество элементов, которые
представляются в виде произведения не более чем $k$
элементов из $X.$

Основная идея работы Дж.~Бергмана \cite{Berg}
заключается в том, что некоторая степень $X^m$ любого порождающего
множества $X$ бесконечной симметрической группы
$\sym \Omega$ содержит стабилизатор вида $\Gamma_{(I),\{J\}}$
(здесь и далее любой символ вида $\Gamma_{*_1,*_2}$ указывает
на пересечение подгрупп $\Gamma_{*_1}$
и $\Gamma_{*_2}$).
где $\Gamma=\sym \Omega,$ а $\Omega = I \cup J$ ---
разбиение $\Omega$ на половинные подмножества.
Далее Дж.~Бергман пользуется тем, что пара стабилизаторов
$\Gamma_{(I),\{J\}} \cup \pi \Gamma_{(I),\{J\}}\pi\inv,$
где $\pi$ --- подходящая перестановка, порождает
$\sym \Omega$ (факт, быстро вытекающий из результатов
работы \cite{DiNeuTho}) и показывает, что
соответствующая ширина равна $3.$
Ясно, что тогда $\sym \Omega \subseteq
X^{3(2k+m)},$ где $k$ --- длина перестановки
$\pi$ по отношению к порождающему множеству
$X.$ Отметим, что бесконечномерные линейные
группы над телами тоже могут быть порождены
парами стабилизаторов вида $\Gamma_{(*),\{*\}}$
\cite{Ev,Macph,To_GL(V)}.

Оказывается, что ключевой результат Дж.~Бергмана может быть
обобщен на порождающие множества групп автоморфизмов
произвольных относительно свободных алгебр бесконечного ранга.
Пусть $F$ --- относительно свободная алгебра
бесконечного ранга. Назовем подалгебру $U$
алгебры $F$ {\it половинной}, если существует
разложение $F=U*V$ в свободное произведение, такое,
что оба свободных множителя $U,V$ изоморфны $F.$
Если теперь $X$ --- это порождающее множество
$\Gamma=\aut F,$ то, как показал
автор в \cite{To_BP}, подходящая степень $X$
содержит стабилизатор вида $\Gamma_{(U),\{V\}},$
где $F=U*V$ --- разложение $F$ в свободное
произведение половинных подалгебр.

Обобщая известные примеры \cite{MN,Berg,Ev,Macph}, вводим следующее
определение \cite{To_BP}. Пусть $\mathfrak V$ --- многообразие
алгебр. Говорим, что $\mathfrak V$ является
{\it BMN-многообразием} (Berg\-man-Macpherson-Neumann
variety), если
для каждой
свободной алгебры $F$ бесконечного ранга,
принадлежащей $\mathfrak V,$ группа автоморфизмов $\Gamma=\aut F$
порождается парой стабилизаторов
\begin{equation}
\Gamma_{(U_1),\{U_2*W\}} \text{ и  }
\Gamma_{(U_2),\{U_1*W\}},
\end{equation}
где $F=U_1 * U_2 * W$ --- разложение $F$
в свободное произведение половинных
подалгебр.

Рассмотренные выше примеры могут быть тогда
интерпретированы в том духе, что многообразие
всех алгебр в пустом языке и любое многообразие
векторных пространств --- это BMN-многообразия.

Используя вышеприведенный результат
о стабилизаторах в степенях порождающих
множеств, нетрудно доказать, что группа автоморфизмов
$\aut F$ свободной алгебры $F$ бесконечного ранга из
BMN-многообразия является группой конечной
ширины в том и только в том случае, если ширина $\aut F$ относительно
любого порождающего множества вида
(\theequation) конечна.  Справедлив также
следующий результат.

\begin{prp}[\cite{To_BP}] \label{NiceProps}
Пусть $F$ --- свободная
алгебра бесконечного ранга из BMN-многообразия.
Тогда группа $\aut F$ порождается инволюциями,
совпадает со своим коммутантом и имеет конфинальность,
строго большую, чем ранг $F.$ Если, кроме того,
группа $\aut F$ является группой конечной ширины, то
группа $\aut F$ имеет свойство FA.
\end{prp}

Напомним, что данная группа $G$ имеет
свойство FA, если всякое действие
$G$ на дереве без инверсий имеет
неподвижную точку \cite[раздел 6.1]{JPS}.

В работе \cite{To_BP} автор показал,
что группа автоморфизмов произвольной
бесконечнопорожденной свободной нильпотентной
(свободной абелевой) группы является
группой конечной ширины. Схема доказательства того, что группа
автоморфизмов $\aut A$ бесконечнопорожденной
свободной абелевой группы $A$ является
группой конечной ширины такова \cite{To_BP}: выбирается
какой-нибудь базис $\cX$ группы $A$,
строятся половинные прямые
слагаемые $U_1,U_2,W,$ порожденные
половинными подмножествами, дающими разбиение
$\cX,$ и проверяется, что
длины всех $\cX$-блочно-унитреугольных
автоморфизмов $A$ по отношению к порождающему
множеству
$$
\Gamma_{(U_1),\{U_2+W\}} \cup
\Gamma_{(U_2),\{U_1+W\}}
$$
ограничены сверху; после этого применяется
\ref{Swan}. Таким образом, группа
$\aut A$ является группой конечной ширины, а многообразие
всех абелевых групп является BMN-многообразием.

Используя результат о свободных абелевых группах,
доказываем индукцией по $c,$ что группа автоморфизмов
$\Gamma =\aut N$ произвольной бесконечнопорожденной
свободной нильпотентной группы ступени $c$ порождается парой
стабилизаторов вида (\theequation) и что
соответствующая ширина конечна. Следовательно, любое
многообразие $\frak N_c$ всех нильпотентных групп
ступени $\le c$ является BMN-мнообразием.

По-видимому, подобные рассуждения трудно
приспособить для изучения групп автоморфизмов
бесконечнопорожденных свободных групп.

\begin{que}
Является ли многообразие всех групп BMN-многообразием$?$
\end{que}

Если ответ на последний вопрос положителен,
то по \ref{NiceProps} группа автоморфизмов
любой бесконечнопорожденной свободной группы
совпадает со своим коммутантом. Выше мы отмечали,
что свободная группа $F$ {\it счетно-бес\-ко\-неч\-ного} ранга
имеет свойство малого индекса, и потому
группа автоморфизмов $F$ совпадает
со своим коммутантом \cite{BrEv,BrRom}. Таким образом,
ответ на следующий вопрос, --- который, конечно, можно
рассматривать вне зависимости от вопроса \theque, ---
положителен для свободных групп счетно-бесконечного
ранга.

\begin{que}
Пусть $F$ --- бесконечнопорожденная свободная
группа. Верно ли, что группа $\aut F$ совпадает
со своим коммутантом$?$
\end{que}

Вопрос Дж.~Бергмана о группах автоморфизмов свободных
групп бесконечного ранга представляется трудным. Тем
не менее, как показал автор в работе \cite{To_BP} группа
автоморфизмов свободной группы счетно-бесконечного
ранга действительно является группой конечной ширины.

Опишем идею ответа на вопрос Дж.~Бергмана
для свободных групп счетно-бесконечного ранга
(всюду в данном абзаце используются определения
и обозначения, введенные на стр. \pageref{0-1Cond}
при обсуждении результатов из работы \cite{Hodges_et_al}).
Пусть $F$ --- свободная группа счетно-бесконечного ранга;
напомним, что согласно результатам из работы \cite{BrEv}
группа $\aut F$ богата генерическими автоморфизмами.
Предполагаем, что найдется порождающее множество
$X$ группы $\aut F,$ по отношению к которому эта
группа имеет бесконечную ширину. Тогда можно, действуя как
и в указанных работах, реализовать игру с генерическими
автоморфизмами группы $F,$ строя индукцией по элементам
$s$ бинарного дерева $2^{< \omega}$ генерические автоморфизмы
$\gamma_s, g_{s\hatin0},g_{s\hatin1},$
удовлетворяющие условиям (i-iii), приведенными
на стр. \pageref{0-1Cond}, и следующему 0-1 условию:
$$
g_{s\hatin0} \in X^n \text{ и }
g_{s\hatin1} \not\in X^{3n}.
$$
Из формулы \eqref{0-1Cond} можно  вывести,
что если функции $\sigma,\pi \in 2^\omega$
таковы, что для некоторого $s \in 2^n$
мы имеем, что $\s$ продолжает $s\hatin0,$
а $\pi$ продолжает $s\hatin1,$ то
длина автоморфизма $\gamma_\sigma\gamma_\pi\inv$
по отношению к $X$ строго больше $n.$ В частности,
существует несчетно много попарно различных
автоморфизмов вида $\gamma_\s,$
где $\s \in 2^\omega,$ а, значит, некоторая
степень $X^n$ порождающего множества $X$
тоже содержит несчетно много автоморфизмов
указанного вида. В таком случае можно
найти $m \ge 2n$ и элементы $\gamma_\s,
\gamma_\pi$ из $X^n,$ для которых имеем:
$$
\sigma \restr m=\tau\restr m \text{ и } \sigma\restr m+1\ne \tau\restr m+1.
$$
Но тогда, как мы замечали выше, длина
автоморфизма $\gamma_\s\gamma_\pi\inv \in X^{2n}$
по отношению к $X$ превосходит $m,$ противоречие.

\begin{que}[Дж.~Бергман]
Верно ли, что группа автоморфизмов любой
бесконечнопорожденной свободной группы
является группой конечной ширины$?$
\end{que}

То обстоятельство, что такие ручные группы
автоморфизмов относительно свободных групп
бесконечного ранга, как группы автоморфизмов
бесконечнопорожденных свободных абелевых и свободных
нильпотентных групп тоже оказались группами конечной
ширины (вне зависимости от ранга), позволяет надеяться
на то, что и в общем случае группы автоморфизмов
бесконечнопорожденных свободных групп будут группами
конечной ширины.

Мы видели ранее, что конфинальность группы автоморфизмов
любой бесконечнопорожденной свободной нильпотентной
(свободной абелевой) группы $N$ строго больше, чем
$|N|.$ Это позволяет надеяться, что результат Р.~Брайана и
Д.~Эванса о несчетной конфинальности группы автоморфизмов
свободной группы счетно-бесконечного ранга можно
обобщить.

\begin{que}
Пусть $F$ --- свободная группа бесконечного
ранга. Превосходит ли конфинальность группы
$\aut F$ ранг группы $F?$
\end{que}

Согласно \ref{NiceProps}, группа автоморфизмов
относительно свободной алгебры $F$ из BMN-многообразия
порождается инволюциями. Более точно, если $\cB$ ---
это какой-нибудь базис $F,$ то группа $\aut F$ совпадает с
нормальным замыканием инволюции $\pi,$ которая
действует на $\cB$ как перестановка порядка два,
такая, что ее носитель равномощен множеству
неподвижных элементов \cite{To_BP}. Так как
многообразие всех абелевых групп и все
многообразия $\frak N_c$ являются
BMN-многообразиями, то группа автоморфизмов бесконечнопорожденной
свободной абелевой (свободной нильпотентной
группы) есть нормальное замыкание некоторой
инволюции (ср. с цитированным ранее результатом
из \cite{BrRom}, утверждающим, что нормальное
замыкание подходящего "перестановочного"\
автоморфизма относительно свободной алгебры
$F$ бесконечного ранга всюду плотно в $\aut F$).

\begin{que}
Пусть $F$ --- бесконечнопорожденная
свободная группа. Верно ли, что
$\aut F$ является нормальным замыканием
подходящей инволюции {\em(}некоторого
автоморфизма группы $F${\em)}$?$
\end{que}

Любопытно отметить, что если $F$ --- неабелева свободная
группа конечного ранга, то $\aut{F}$
есть нормальное замыкание инволюции, которая
действует на некотором базисе $F$ как перестановка
с единственным циклом длины два \cite{BriVo_Homs}.
М.~Брайдсон и К.~Фогтманн, авторы последнего
результата, применили его для описания возможных
гомоморфных образов группы $\aut{F}.$

Любой гомоморфный образ группы конечной ширины
--- это тоже группа конечной ширины. Поэтому,
скажем, любой гомоморфный образ группы автоморфизмов
свободной группы $F$ счетно-бесконечного ранга обязан
быть группой конечной ширины, что, как и в случае
свободных групп конечного ранга, накладывает серьезные
ограничения на гомоморфные образы группы $\aut F.$
С этой точки зрения  две следующие проблемы
имеют важное значение.

\begin{que}[Дж.~Бергман, \cite{Berg}]
Существует ли счетно-бесконечная группа
конечной ширины$?$
\end{que}

\begin{que}[М.~Дросте, Р.~Гёбель, \cite{DrGo}]
Существует ли группа конечной ширины,
имеющая счетную конфинальность$?$
\end{que}

Напомним, что если $G$ --- это относительно
свободная группа, то $\ia G$ --- это группа
так называемых IA-автоморфизмов, т.е. автоморфизмов,
действующих тождественно на абелизации $G.$
Хорошо известно, что группа $\ia{F}$ всех IA-автоморфимов
свободной группы $F$ конечного ранга $n$ является
нормальным замыканием (в группе $\aut{F})$
следующего $\iat$-ав\-то\-мор\-физма:
$$
\begin{cases}
\a x_1 = x_1, \\
\a x_2 = x_1 x_2 x_1\inv, \\
\a x_i = x_i, \quad i \ne 1,2,
\end{cases}
$$
где $x_1,\ldots,x_n$ --- произвольный базис
$F.$

Было бы интересно проверить, справедлив ли аналогичный
результат для бесконечнопорожденных свободных групп.

\begin{que}
Пусть $F$ --- бесконечнопорожденная свободная
группа. Существует ли $\iat$-автоморфизм группы
$F$ {\em(}некоторое "небольшое"\ множество
$\iat$-автоморфизмов{\em)}, нормальное замыкание
которого в группе $\aut F$ есть подгруппа $\ia F?$
\end{que}

Пусть $F$ --- произвольная неабелева свободная группа.
Укажем пример множества $\iat$-автоморфизмов,
нормальное замыкание которого совпадает с группой $\ia F,$ в
случае, если ранг $F$ конечен, и нормальное замыкание
которого совпадает, возможно, с группой $\ia F$ и в общем
случае. Назовем {\it симметрией} инволюцию группы $F,$
которая обращает все элементы некоторого базиса $F.$
Ясно, что произведение двух симметрий индуцирует
тривиальный автоморфизм абелизации $F.$ Тогда, если
ранг $F$ конечен, то нормальное замыкание множества всех
произведений пар симметрий в группе $\aut F$
есть подгруппа $\ia F$, ибо указанный
выше автоморфизм $\a$ легко представить в виде
произведения двух симметрий.

\section{Автоморфизмы}

Классификация изоморфизмов (автоморфизмов) для
основных типов линейных групп над телами,
полученная в классических работах Ж.~Дьедонне,
К.~Риккарта и других авторов, обусловила появление
важной работы Л.-К.~Хуа и И.~Райнера \cite{HuaRei} 1951 года,
в которой было найдено описание автоморфизмов
групп автоморфизмов свободных абелевых групп
конечного ранга (унимодулярных групп $\text{GL}(n,\Z)$). Группы автоморфизмов свободных
абелевых групп  занимают особое место среди групп
автоморфизмов относительно свободных групп,
ибо эти группы могут рассматриваться как группы
автоморфизмов свободных модулей (над $\Z$), что
позволило Л.-К.~Хуа и И.~Райнеру применить ряд методов,
разработанных для изучения групп автоморфизмов
векторных пространств над телами.

Пусть $A$ --- свободная абелева групп конечного ранга
$n \ge 3.$ Как показано в работе \cite{HuaRei}, в группе
$\aut{\aut A}$ имеется не более четырех смежных классов
по подгруппе $\inn{\aut A}$ внутренних автоморфизмов
группы $\aut A$: а именно, ровно четыре класса,
если $n$ --- четно, и ровно два класса, если $n$ нечетно.
Как и в случае, скажем, общих линейных групп над
телами, один из нетривиальных классов определяется
(при любом $n$) контраградиентным автоморфизмом $R_1$,
а  два других (только при четном $n$) ---
автоморфизмом $R_2(\sigma)= \det \s \cdot \sigma$ и
автоморфизмом $R_1 R_2.$

Новые результаты об автоморфизмах групп автоморфизмов
относительно свободных групп появились только примерно
через четверть века после выхода работы Л.-К.~Хуа и И.~Райнера.
Появлению этих результатов способствовал ряд гипотез
Г.~Баумслага о башнях автоморфизмов групп,
предложенных в начале семидесятых. В частности, он
сформулировал гипотезу о том, что башня автоморфизмов
свободной группы конечного ранга должна быть очень
короткой и, пожалуй, наиболее известную (и до сих пор
ни подтвержденную, ни опровергнутую)  гипотезу о том,
что башня автоморфизмов всякой нильпотентной группы
без кручения должна обрываться после конечного числа
шагов \cite[проблема 4.9]{Kou}. (Заметим, что для
простоты изложения мы допускаем, как это
принято в последнее время, построение башни
автоморфизмов над любой группой, а не только
над группой без центра, как того требует классическое
определение.)

В серии работ \cite{DFo_abs,DFo_2-step,DFo_CompAutGrs,
DFo_CharSubgrs}, написанных в середине семидесятых,
Дж.~Дайер и Э.~Форманек подтвердили некоторые из гипотез
Г.~Баумслага и описали башни автоморфизмов для достаточно
большого класса относительно свободных групп.
Выяснилось, к примеру, что башня автоморфизмов
неабелевой свободной группы $F$ конечного ранга
является настолько короткой, насколько это вообще
возможно, поскольку как показали Дж.~Дайер и Э.~Форманек в
работе \cite{DFo_abs} имеет место изоморфизм групп $\aut F$ и
$\aut{\aut F}.$ Иными словами, группа $\aut F$
оказалось совершенной (напомним, что группа $G$
называется {\it совершенной}, если центр $G$
\index{совершенная группа} тривиален, и все ее
автоморфизмы --- внутренние). Дж.~Дайер и Э.~Форманек
использовали следующий критерий совершенности
группы автоморфизмов $\aut G$ группы без центра
$G$, принадлежащий У.~Бернсайду:
группа $\aut G$ совершенна в том и только
в том случае, если подгруппа внутренних
автоморфизмов $\inn G$ является характеристической
подгруппой группы $\aut G.$ Отметим, что в работе \cite{DFo_abs}
существенным образом использовался приведенный
выше результат Л.-К.~Хуа и И.~Райнера.

Изучая частный случай вышеприведенной
гипотезы Г.~Баумслага о башнях автоморфизмов
нильпотентных групп без кручения, Дж.~Дайер и Э.~Форманек
установили в работе \cite{DFo_2-step}, что группа автоморфизмов
свободной нильпотентной группы ступени
два, имеющей конечный ранг $r \ge 2,$ является
совершенной при $r \ne 3,$ а при $r=3$ соответствующая башня автоморфизмов
имеет высоту $2.$

В работе \cite{DFo_CharSubgrs} Дж.~Дайер и Э.~Форманек
изучали группы автоморфизмов относительно свободных групп
вида $F/C,$ где $C$ --- характеристическая
подгруппа свободной группы $F$ конечного ранга.
Один из наиболее общих результатов \cite{DFo_CharSubgrs}
говорит о том, что если фактор-группа $F/R,$ где
$R$ --- характеристическая подгруппа $F,$ аппроксимируется
нильпотентными группами без кручения, то
группа автоморфизмов фактор-группы $F/R',$ где $R'$
--- коммутант группы $R,$ является совершенной.
Напомним, что свободные группы аппроксимируются
нильпотентными группами без кручения, и поэтому
последний результат обобщает, на самом деле,
результат о совершенности группы $\aut F.$
Среди относительно свободных
групп аппроксимируемых нильпотентными группами
без кручения можно упомянуть еще, к примеру,
разрешимые группы; поэтому из указанного выше
результата вытекает, что группа автоморфизмов
каждой неабелевой свободной разрешимой группы
является совершенной.

В работе \cite{D_GL} Дж.~Дайер изучала
башни автоморфизмов групп
$\text{SL}(n,\Z)$ и $\text{GL}(n,\Z)$.
Было доказано, что группа автоморфизмов
группы $\text{SL}(n,\Z)$ совершенна
при любом $n \ge 3$, а при $n=2$
башня автоморфизмов $\text{SL}(n,\Z)$ имеет высоту
$8.$ Группа автоморфизмов
группы $\text{GL}(n,\Z)$ также совершенна
при любом нечетном $n \ge 3.$ Если $n \ge 4$
четно, то высота
башни автоморфизмов группы $\text{GL}(n,\Z)$
равна 5, а при $n=2$ --- равна $4.$

В 1991 Э.~Форманек \cite{Fo} усилил результат
из \cite{DFo_abs}, явным образом указав на причину
того, что подгруппа $\inn F$ является характеристической
подгруппой группы автоморфизмов $\aut F$ неабелевой
свободной группы $F$ конечного ранга: ключевой
результат из \cite{Fo} гласит, что $\inn F$ --- это
единственная нормальная свободная подгруппа группы $\aut F,$
имеющая ранг, совпадающий с рангом $F.$ Еще одно
доказательство совершенности $\aut F$, найденное в том
же году Д.~Г.~Храмцовым \cite{Khr2}, основывалось на
полученной им классификации конечных групп,
реализуемых как подгруппы групп автоморфизмов
конечнопорожденных свободных групп \cite{Khr_MZ,Khr1}.
Кроме того, Д.~Г.~Храмцов показал, что и группа
$\out F$ внешних автоморфизмов группы $F$ является
совершенной при условии, что $\rank F \ge 3$ (напомним,
что группа внешних автоморфизмов свободной группы
ранга два изоморфна группе $\text{GL}(2,\Z)),$
и, в частности, не является совершенной).

М.~Брайдсон и К.~Фогтманн передоказали в 2000 году
\cite{BriVo} результат Д.~Г.~Храм\-цо\-ва о
совершенности $\out F,$ где $F$ --- неабелева
свободная группа ранга не меньше $3$ и предложили еще
одно (четвертое по счету) доказательство совершенности
группы автоморфизмов неабелевой свободной группы
конечного ранга. М.~Брайдсон и К.~Фогтманн
использовали действие группы $\out F$ на так
называемом внешнем пространстве (Outer Space),
специальном и очень полезном для приложений
комбинаторном объекте, введенном К.~Фогтманн и
М.~Каллером в \cite{CuVo}.

В своей диссертации \cite{Kass} (2003 год)
М.~Кассабов показал, что башня автоморфизмов свободной
нильпотентной группы конечного ранга обрывается после
конечного числа шагов, решив, тем самым, для этих
групп проблему Г.~Баумслага. Например, если
$N$ --- свободная нильпотентная группа конечного
ранга, такая, что группа $\aut N$ имеет тривиальный
центр, то башня автоморфизмов над $N$ имеет
высоту $\le 3.$ Основная идея работы М.~Кассабова заключается во
вложении каждого этажа соответствующей башни
автоморфизмов в качестве решетки в подходящую группу
Ли. Это позволило М.~Кассабову свести задачу об
описании башни автоморфизмов свободной нильпотентной
группы к задаче об описании башни производных
некоторой свободной нильпотентной алгебры Ли.
Представляется вероятным, что
некоторые результаты М.~Кассабова можно усилить.

Результаты М.~Кассабова были стимулированы важной
работой Э.~Форманека \cite{Fo_Centers} 2002 года, описавшим центры групп
автоморфизмов свободных нильпотентных групп конечного
ранга. Пусть $N_{r,c}$ --- свободная
нильпотентная группа конечного ранга $r$ и ступени
нильпотентности $c.$ Заметим, что проблема описания
центров групп $\aut{N_{r,c}}$ связана со следующим интересным вопросом
А.~Г.~Мясникова: для каких $r,c$ в группе $N_{r,c}$
найдутся нетривиальные элементы, фиксируемые всеми
автоморфизмами $N_{r,c}?$ (см. \cite[проблема
N1]{BaMyaSh}). Как показал Э.~Форманек в
работе \cite{Fo_Centers}, группа $N_{r,c}$ обладает такими
элементами тогда и только тогда, когда (a) $r=2$ или
$r=3$ и $c=2kr$, $k \ge 2$ и (b) $r\ge 4$ и $c=2kr$,
$k\ge 1$. В.~В.~Блудов построил в явном виде примеры
нетривиальных элементов, фиксируемых всеми
автоморфизмами групп $N_{2,4k},$ где $k \ge 2$
\cite{BaMyaSh}.  Группа ${\rm Aut}(N_{r,c})$ обладает
нетривиальным центром тогда и только тогда, когда
$c=2kr+1$ и $k\ge 1$.

Все цитированные выше статьи об автоморфизмах групп
автоморфизмов относительно свободных групп конечного
ранга так или иначе существенным образом используют
конечность ранга: например, при использовании действия
автоморфизмов на порождающих, как в
\cite{DFo_abs,DFo_2-step,Khr2,BriVo} или при
использовании результата Л.-К.~Хуа и И.~Райнера, существенно
использующего матричную технику, как в
\cite{DFo_abs,DFo_2-step,DFo_CompAutGrs,
DFo_CharSubgrs,Fo}. Представляется, однако,
целесообразным обобщение этих результатов на группы
автоморфизмов относительно свободных групп
бесконечного ранга. В цикле работ
\cite{ToJLM1,ToCamb,ToContMat1,ToContMat2} автор начал
реализацию этой программы.

Основной результат работы \cite{ToJLM1} обобщает
результаты из \cite{DFo_abs}:  оказывается, что группа
автоморфизмов $\aut F$ неабелевой свободной группы $F$
является совершенной вне зависимости от того, является
ли ранг $F$ конечным или бесконечным.  Это, как и в
случае конечного ранга следует из того, что подгруппа
$\inn F$ является характеристической подгруппой группы
$\aut F,$ что, в свою очередь, обусловлено следующим
результатом.

\begin{prp} \label{Def-of-Conjs-by-PPE}
Пусть $F$ --- неабелева свободная
группа. Тогда семейство всех сопряжений
{\em(}внутренних автоморфизмов $F${\em)} посредством  степеней
примитивных элементов $F$ является определимым
без параметров в логике первого порядка
подмножеством группы $\aut F.$
\end{prp}

Напомним, что подмножество $D$ структуры $\cM$
{\it определимо без параметров} в $\cM,$ если найдется
формула первого порядка $\chi$ в языке $\cM,$ такая,
что $D=\chi(\cM)$; ясно, что если $\cM$ --- это
группа, то всякое определимое без параметров
подмножество $\cM$ порождает характеристическую
подгруппу.

Доказательство теоремы \ref{Def-of-Conjs-by-PPE}
опирается на описание автоморфизмов простого
порядка, найденное Дж.~Дайер и П.~Скоттом
в работе \cite{DSc}. А именно, используя
результаты из \cite{DSc}, можно показать,
что класс сопряженности инволюций, называемых
в \cite{ToJLM1} квази-сопряжениями, является
определимым без параметров классом сопряженности
группы $\aut F,$ а потом, выбрав некоторое
квази-сопряжение $\varphi,$ построить $\varphi$-определимое
подмножество $\Pi(\varphi),$ централизатор которого таков, что
его элементы, не являющиеся инволюциями,
суть сопряжения всеми степенями некоторого
примитивного элемента $F.$

Воспроизведем теперь детали набросанной
выше схемы доказательства основного
результата из работы \cite{ToJLM1}. Пусть $x$ --- примитивный
элемент группы $F$ и $C$ --- свободный множитель
$F,$ такой, что $F=\str x * C.$ Тогда инволюция
$\f \in \aut F,$ обращающая $x,$ и действующая
на $C$ как сопряжение посредством $x,$ называется
{\it квази-сопряжением.} Назовем подмножество
$S$ данной группы $G$ {\it анти-коммутативным},
если все его элементы попарно неперестановочны.
В случае, если $\rank F > 2$ класс сопряженности
квази-сопряжений --- это
единственный антикоммутативный класс сопряженности
инволюций, такой, что для всякого другого
антикоммутативного класса сопряженности инволюций
$K'$, все инволюции в $KK'$ попарно
сопряжены. В случае, если $\rank F=2$ класс сопряженности
квази-сопряжений является
единственным антикоммутативным классом
сопряженности инволюций, элементы которого
не являются квадратами. Если $\varphi$ --- это
квази-сопряжение, то указанное выше множество
$\Pi(\varphi)$ состоит из всех автоморфизмов
вида $\sigma\sigma',$ где $\sigma,\sigma'$
--- элементы из централизатора
$\f,$ такие, что $\sigma$ и $\sigma'$ сопряжены
в группе $\aut F.$

Существуют свидетельства в пользу того, что
ответ на следующий вопрос положителен.

\begin{que}
Пусть $F$ --- свободная группа бесконечного
ранга. Совершенна ли группа $\out F?$
\end{que}

В работе \cite{ToCamb} автор обобщил результат
Дж.~Дайер и Э.~Форманека из работы \cite{DFo_2-step}.
Основной результат работы \cite{ToCamb} гласит,
что группа автоморфизмов бесконечнопорожденной
свободной нильпотентной группы $N$ ступени два совершенна.
Доказательство, как и в \cite{DFo_2-step}, начинается
с доказательства характеристичности подгруппы $\ia N$
в группе $\aut N.$ На самом деле, если $G$ --- любая
свободная нильпотентная группа ступени два, то
подгруппа $\ia G$ является определимой без параметров
подгруппой $\aut G$ в логике первого порядка (ср. с
аналогичным результатом из \cite{DFo_2-step}, утверждающим, что
если $G$ имеет конечный ранг, то $\ia G$ --- это радикал Хирша-Плоткина группы $\aut
G$). На следующем этапе доказывается определимость
в логике первого порядка подгруппы $\inn N$
в группе $\aut N.$ Результат У.~Бернсайда,
однако,  применить нельзя, ибо $N$ имеет нетривиальный
центр; возникающую трудность можно, тем не менее,
обойти, показав, что в группе $\aut N$
можно восстановить теоретико-групповыми
методами семейство всех примитивных элементов
группы $N.$

\begin{que}
Описать автоморфизмы групп автоморфизмов
свободных абелевых {\em(}нильпотентных, разрешимых{\em)}
групп бесконечного ранга.
\end{que}

В работах \cite{ToContMat1,ToContMat2}, посвященных
решению одной логической проблемы, которая
будет обсуждаться в следующем параграфе, получены
результаты, позволяющие надеяться, что все
автоморфизмы группы автоморфизмов свободной
абелевой группы бесконечного ранга являются
внутренними, и что группа автоморфизмов
любой свободной нильпотентной группы бесконечного
ранга является совершенной.

\section[Выразительная сила теорий первого порядка]{Выразительная сила теорий первого порядка и классификация
элементарных типов}

В начале семидесятых широкое внимание логиков привлек
вопрос Дж.~Ис\-бел\-ла о классификации элементарных типов
бесконечных симметрических групп (классификации этих
групп с точностью до элементарной эквивалентности).
Обобщая результаты, полученные рядом авторов, С.~Шелах
дал в работах \cite{Sh1,Sh2} окончательное решение
этой проблемы. Для пояснения формулировок результатов
из \cite{Sh1} нам потребуется следующее определение.

Пусть $\{T^0_i : i \in {\bold I}\}$ и  $\{T^1_i : i \in{\bold I}\}$
--- семейства теорий в логиках
$\cL_0$ и $\cL_1,$ соответственно.  Говорят,
что теория $T^0_i$ {\it синтаксически интерпретируется} в
теории $T^1_i$ {\it равномерно {\em(}единообразно{\em)} по} $i \in {\bold I},$
если существует отображение
$\,{}^*\,$ множества всех $\cL_0$-предложений во
множество всех $\cL_1$-предложений, такое, что
для каждого $\cL_0$-предложения $\chi$ и для каждого $i \in {\bold
I},$ $\chi \in T^{0}_i,$ если и только если $\chi^* \in T^1_i$ \cite{BSh,Ersh,Ho}.
Если, дополнительно,
теория $T^1_i$ синтаксически интерпретируется
в теории $T^0_i$ равномерно по $i \in {\bold I},$ то
теории $T^0_i, T^1_i$ называют {\it
взаимно синтаксически интерпретируемыми равномерно по $i \in
{\bold I}.$} Если $T^0_i=\Theo(\cM_i,\cL_0)$ и
$T^1_i= \Theo(\cN_i,\cL_1),$ т.е. если рассматриваемые
теории являются теориями некоторых структур
в логиках $\cL_0$ и $\cL_1,$ соответственно,
то естественным достаточным условием
для синтаксической интерпретируемость теории
$T_i^0$ в теории $T_i^1$ является единообразная интерпретируемость
(интерпретируемость с $\cL_1$-определимыми параметрами)
структуры $\cM_i$ в структуре $\cN_i$ средствами
логики $\cL_1$ для всех $i \in \bold I$ \cite{BSh,Ersh,Ho}. Кроме того,
в данном случае имеем, что $\cN_i \equiv_{\cL_1}
\cN_j$ влечет $\cM_i \equiv_{\cL_0}
\cM_j$ для всех $i,j \in \bold I.$

С.~Шелах показал в \cite{Sh1}, что элементарная теория
симметрической группы $\sym{\aleph_\alpha}$ над
кардиналом $\aleph_\alpha$ взаимно синтаксически
интерпретируема с теорией двухсортной структуры
$\str{\alpha,\lambda_\alpha;<}$ в логике $\bold
L_2((2^{\aleph_0})^+),$ где $\lambda_\alpha$ --- кардинал
$\min(\aleph_\alpha,2^{\aleph_0}),$ рассматриваемый как
множество без структуры, а $<$ --- отношение полного
порядка ординала $\alpha.$ Логика $\bold L_2(\varkappa^+)$
(C.~Шелах), где $\varkappa$ --- некоторый кардинал,
 --- это фрагмент полной логики второго порядка,
допускающий квантификацию по отношениям универсума
мощности не выше $\varkappa.$ Результат из \cite{Sh1}
резко контрастирует с известным результатом М.~Рабина,
утверждающим что элементарная теория полугруппы $\text{End}(\aleph_\alpha)$ всех
отображений кардинала
$\aleph_\alpha$ в себя взаимно синтаксически
интерпретируема с теорией $\Theo_2(\aleph_\alpha)$ --- теорией
 кардинала $\aleph_\alpha,$ рассматриваемого как
множество без структуры, в полной логике второго порядка $\bold L_2$.
Неформально, сравнивая результаты М.~Рабина и С.~Шелаха,
можно сказать, что выразительная сила элементарной
теории полугруппы $\text{End}(\aleph_\alpha)$
превосходит выразительную силу элементарной теории
группы $\text{Sym}(\aleph_\alpha).$

Работа \cite{Sh1} послужила для С.~Шелаха отправной
точкой для ряда важных работ. Одна из них --- работа \cite{Sh3}
1977 года может рассматриваться как значительное обобщение
упомянутого выше результата М.~Рабина --- на полугруппы
эндоморфизмов свободных объектов в многообразиях
алгебр. Пусть  $F_\varkappa=F_\varkappa(\mathfrak V)$ --- свободная
алгебра бесконечного ранга $\varkappa$
из многообразия алгебр $\mathfrak V$ в языке $L.$ Тогда
если $\varkappa > |L|,$ то элементарная
теория полугруппы эндоморфизмов алгебры $F_\varkappa$
синтаксически интерпретирует теорию кардинала
$\varkappa,$  рассматриваемого как множество
без структуры, в полной логике второго
порядка \cite{Sh3}.

С.~Шелах замечает в \cite{Sh3}, что естественно
изучать вопрос о том, когда выразительная
сила элементарной теории {\it группы автоморфизмов} свободной
алгебры $F_\varkappa$ сравнима с выразительной
силой элементарной теории полугруппы эндоморфизмов $F_\varkappa.$
Через более чем двадцать лет С.~Шелах вновь вернулся
к этому вопросу, включив его в список проблем
в обзоре \cite[проблема 3.14]{ShMathJap}, и предложив
описать многообразия алгебр $\mathfrak V,$
группы автоморфизмов свободных алгебр $F_\varkappa(\mathfrak V)$ которых
интерпретируют средствами логики первого порядка
теорию кардинала $\varkappa$ в полной логике
второго порядка для всех (или, возможно,
всех достаточно "больших") бесконечных кардиналов
$\varkappa.$ Конечно, ситуация с группами автоморфизмов
свободных алгебр выглядит гораздо более сложной, ибо,
несмотря на то, что полугруппы автоморфизмов свободных алгебр
могут быть исключительно сложными, они, все же, если
предельно обобщить существо результатов С.~Шелаха из \cite{Sh3}, являются "комбинаторными"\ объектами.
Насколько известно автору, никаких более или менее общих
результатов по проблеме С.~Шелаха нет. В \cite{ShMathJap}
С.~Шелах предлагает схему изучения проблемы для
класса многообразий, который он называет классом
Aut-разложимых многообразий. Предполагается, что группы автоморфизмов
свободных алгебр бесконечного ранга из этих многообразий ведут себя во многом
так же, как и бесконечные симметрические группы (см. подробности
в \cite{ShMathJap}).

В своей кандидатской диссертации \cite{Tolsty'sPh},
результаты из которой опубликованы в работе
\cite{ToAPAL}, автор нашел
решение проблемы С.~Шелаха для многообразий векторных
пространств над телами. Пусть $V$ ---
векторное пространство бесконечной размерности
$\varkappa$ над телом $D.$ Тогда если $\varkappa > |D|,$ то
элементарная теория группы $\gl V$ взаимно синтаксически
интерпретируема с теорией двухсортной структуры
$\str{\varkappa,D},$ основными отношениями которой
являются только основные отношения тела $D,$ в
логике $\bold L_2(\varkappa^+).$ Заметим, что
если $D$ --- поле, то указанный результат верен
для любых бесконечных кардиналов $\varkappa.$

В работах \cite{ToJLM2,ToContMat1,ToContMat2}
автор рассматривал проблему С.~Шелаха для классических
многообразий групп. Результаты из упомянутых
работ можно суммировать следующим образом. Пусть
$\mathfrak V$ --- это либо многообразие всех групп,
либо многообразие всех абелевых групп, либо любое
многообразие $\frak N_c$ всех нильпотентных групп ступени
$\le c.$
Тогда если $F_\varkappa(\mathfrak V)$ --- свободная
группа $\mathfrak V$ бесконечного ранга $\varkappa,$
то элементарная теория группы $\aut{F_\varkappa(\mathfrak V)}$
и теория в полной логике второго порядка кардинала $\varkappa$
как множества без структуры взаимно синтаксически
интерпретируемы, равномерно по $\varkappa.$
Таким образом, для указанных многообразий групп
проблема С.~Шелаха решается положительно. В качестве
следствия получаем классификацию элементарных типов
групп автоморфизмов:
$$
\aut{F_\varkappa(\mathfrak V)} \equiv
\aut{F_\lambda(\mathfrak V)} \iff
\varkappa \equiv_{\bold L_2} \lambda.
$$

В работе \cite{ToJLM2} рассматривается случай
многообразия всех групп. Пусть $F$ ---
бесконечнопорожденная свободная группа.
В силу того, что ранг
$F$ бесконечен, любой элемент $F$ есть, как легко
видеть, произведение двух примитивных элементов.
Тогда из \ref{Def-of-Conjs-by-PPE} легко вытекает, что
подгруппа $\inn F$ всех сопряжений является
определимой без параметров в логике первого порядка
подгруппой группы $\aut F.$ Таким образом, получаем,
что в группе $\aut F$ интерпретируется двухсортная
структура $\str{\aut F,F},$ основные отношения
которой --- это групповые операции на $\aut F$
и $F,$ а также предикат, задающий действие $\aut F$ на
$F.$

Несколько отступая от темы, заметим, что было бы очень
интересно проверить, справедлив ли аналогичный
результат для свободных групп конечного ранга.

\begin{que}
Пусть $F_n$ --- неабелева свободная группа
конечного ранга $n.$ Является ли подгруппа
$\inn{F_n}$ определимой подгруппой группы
$\aut{F_n}?$
\end{que}

Из результатов работы \cite{BST}, написанной автором
совместно с В.~Г.~Бардаковым и В.~Э.~Шпильрайном,
вытекает, что, по крайней мере, опираясь только на
определимость сопряжений степенями примитивных
элементов, ответ на последний вопрос (напрямую)
получить нельзя.  Действительно, один из основных
результатов \cite{BST} гласит, что примитивная ширина
неабелевой свободной группы $F_n,$ т.е. ширина $F_n$
относительно множества всех примитивных элементов
$F_n,$ --- является бесконечной. Больше того, можно
показать, что если $\{x,y\}$ --- это некоторый базис
свободной группы $F_2$ ранга $2,$ то множество
всех сопряжений палиндромами от букв $x,y$ (т.е.
неприводимыми словами в буквах $x,y,$ которые
читаются одинаково как справа налево, так
и слева направо) является определимым с параметрами из определимого
множества подмножеством $\aut{F_2}.$ Однако,
палиндромная ширина любой неабелевой свободной группы
является бесконечной \cite{BST}.

Следующий шаг в решении проблемы С.~Шелаха
для многообразия всех групп --- это интерпретация в структуре
$\str{\aut F,F}$ трехсортной структуры
$\str{\aut{F},F,S},$ где $S$ обозначает множество всех
свободных множителей $F.$ Основные отношения последней
структуры включают в себя: основные отношения структур
$\aut F$ и $F$; отношения, задающие действие группы
$\aut F$ на $F$ и $S;$ отношение принадлежности на
множестве $F \cup S$ и тернарное отношение $R(A,B,C)$
на $S,$ выполняющееся если и только если $A=B*C.$ Интерпретация
существенным образом использует результаты об
инволюциях группы $\aut F,$ полученные в
\cite{ToJLM1}.

Затем в структуре
$\str{\aut{F},F,S}$ интерпретируется действие группы на некотором
базисе $F$ или, более формально, в
указанной структуре интерпретируется структура  $\str{\aut F,F,\cB}$ (с
естественными отношениями), где $\cB$ ---
какой-нибудь базис $F.$ Далее, используя достаточно
стандартные методы, в структуре
$\str{\aut F,F,\cB}$ можно интерпретировать структуру
$\str{|F|^{|F|},|F|},$ где $|F|^{|F|}$ --- множество
всех отображений кардинала $|F|$ в себя.  После этого
применяется цитированный выше результат М.~Рабина,
позволяющий интерпретировать в элементарной теории
последней структуры теорию множества $|F|$ в полной
логике второго порядка.  Синтаксическая интерпретация
элементарной теории группы $\aut F$ в теории
$\Theo_2(|F|)$ строится довольно просто.

В работе \cite{ToContMat2}
проблема С.~Шелаха изучается для многообразия $\frak A$
всех абелевых групп. Как было сказано
выше, для этого многообразия проблема
С.~Шелаха тоже решается положительно.
Ясно тогда, что и для любого другого многообразия
групп $\mathfrak V$ проблема С.~Шелаха решается положительно,
если абелизации свободных групп из $\mathfrak V$ являются
свободными абелевыми группами, и если в группе
автоморфизмов $\aut{F(\mathfrak V)},$ где $F(\mathfrak V)$ ---
бесконечнопорожденная свободная группа из $\mathfrak V,$
можно средствами логики первого порядка
интерпретировать группу автоморфизмов абелизации
$F(\frak V).$ Подобный подход, как показывается
в статье автора \cite{ToContMat2}, может быть применен
к любому многообразию $\frak N_c$ всех нильпотентных
групп ступени $\le c.$

Пусть $A$ --- бесконечнопорожденная свободная абелева
группа. Обсудим схему интерпретации теории
множества $|A|$ в полной логике второго порядка  в
элементарной теории группы $\aut A$ из
\cite{ToContMat2}.

Сначала в $\aut A$ восстанавливаются средствами логики
первого порядка фрагменты геометрии $\Z$-модуля $A$: а
именно, средствами
логики первого порядка строится интерпретация в $\aut A$  семейства $\cD^1(A),$
состоящего из всех прямых слагаемых группы $A,$
имеющих ранг или коранг $1.$ Для этой цели,
как и в идейно близкой к рассматриваемой работе
работе \cite{ToAPAL},
используются классические методы из теории
линейных групп, основанные на геометрических свойствах
инволюций.

На следующем этапе, используя действие $\aut A$ на
$\cD^1(A),$ можно показать, что главная
конгруэнц-подгруппа $\Gamma_2(A)$ уровня два,
состоящая из всех автоморфизмов $A,$ действующих
тривиально на фактор-группе $A/2A,$ является
определимой без параметров подгруппой группы $\aut A.$
Фактор-группа $\aut A/\Gamma_2(A)$ изоморфна общей
линейной группе векторного пространства $A/2A$ над
полем $\Z_2.$ Таким образом, группа $\aut A$
интерпретирует средствами первого порядка группу
$\gl{|A|,\Z_2}=\aut{A/2A}.$ Элементарная теория последней группы,
как мы замечали выше, синтаксически интерпретирует
теорию $\Theo_2(|A|).$

Как уже говорилось, естественная схема решения проблемы С.~Шелаха
для многообразий $\frak N_c$  заключается в интерпретации
в группе $\aut N,$ где $N$ --- бесконечнопорожденная
свободная нильпотентная группа, группы автоморфизмов
свободной абелевой группы, имеющей тот же ранг,
что и $N.$ То, что такая схема действительно реализуема,
следует из результатов еще одной работы автора \cite{ToContMat1}.

Инволюции относительно свободных групп, обращающие
элементы некоторого базиса, автор называет в своих
работах {\it симметриями}. В \cite{ToContMat1} показывается,
что семейство всех инволюций, сравнимых по модулю
подгруппы $\ia N$ с симметриями из $\aut N$ является
определимым без параметров семейством $\aut N.$

Отсюда можно довольно легко вывести такой результат:
группа $\aut N$ интерпретирует средствами первого
порядка группу $\aut{N/\gamma_c(N)},$ где
$c$ --- это ступень нильпотентности группы $N,$ а $\gamma_m(N)$
обозначает, как обычно, член нижнего центрального ряда $N$
с номером $m.$ Дело в том, что ядро $K$
гомоморфизма $\aut N \to \aut{N/\gamma_c(N)},$
индуцированного естественным гомоморфизмом
$N \to {N/\gamma_c(N)},$ является определимой
подгруппой $\aut N$:
$$
K = T^+(N) \cup T^-(N),
$$
где $T^+(N)$ (соотв. $T^-(N)$) обозначает множество
всех автоморфизмов $N,$ сохраняемых (соотв. обращаемых)
при присоединенном действии всеми инволюциями,
сравнимыми с симметриями по модулю подгруппы
$\ia N.$

Ясно, что группа $N/\gamma_c(N)$ является
нильпотентной группой ступени $c-1$ и
ее ранг совпадает с рангом группы
$N.$ Теперь простое рассуждение,
использующее индукцию по $c,$ показывает,
что в $\aut N$ интерпретируется
группа $\aut A,$ где $A$ --- свободная
абелева группа, имеющая тот же ранг,
что и $N.$

Как видим, решение проблемы С.~Шелаха для многообразий
абелевых групп и многообразий $\frak N_c$ положительно;
с другой стороны, группы автоморфизмов свободных
групп бесконечного ранга из этих многообразий --- ручные.
Естественен, поэтому, следующий вопрос.

\begin{que}
Пусть $\mathfrak V$ --- многообразие групп, такое,
что группы автоморфизмов свободных групп бесконечного
ранга из этого многообразия --- ручные. Верно ли,
что для многообразия $\mathfrak V$ проблема С.~Шелаха
решается положительно$?$
\end{que}

Как отмечалось в предыдущем параграфе, результаты из
\cite{ToContMat1,ToContMat2} могут быть, по-видимому, применены
для подтверждения гипотез автора об автоморфизмах
групп автоморфизмов свободных абелевых  и свободных
нильпотентных групп, сформулированных в конце
предыдущего параграфа.

\newpage

УДК 512.543.25, 510.67


\bigskip

\bigskip

{Группы автоморфизмов относительно свободных групп бесконечного ранга}

\bigskip

В.~А.~Толстых (Кемеровский государственный университет)

\bigskip

Обозорная статья, посвященная алгебраическим и логическим
результатам  о группах автоморфизмов относительно свободных
групп бесконечного ранга. Обсуждаются следующие вопросы:
топологические методы, результаты о порождающих
множествах, описание автоморфизмов, выразительная
сила элементарных теорий и классификация элементарных
типов. Библиогр. 56 назв.

\end{document}